\documentclass{amsart}
\usepackage{graphicx}

\input{amssym.def}
\input{amssym.tex}

\newcommand{\ben}{\begin{enumerate}}
\newcommand{\een}{\end{enumerate}}

\newcommand{\be}{\begin{enumerate}}
\newcommand{\ee}{\end{enumerate}}
\newcommand{\bq}{\begin{eqnarray*}}
\newcommand{\eq}{\end{eqnarray*}}

\newcommand{\disp}{\displaystyle}
\begin{document}

\title{\bf ON THE NUMBER OF FUZZY SUBGROUPS OF A SYMMETRIC GROUP $ S_5{} $ }
\author{Ogiugo, M.E. and EniOluwafe, M.}
\date{\today}

\maketitle
         Department of Mathematics, University of Ibadan, Ibadan.

 e-mail: ekpenogiugo@gmail.com, michael.enioluwafe@gmail.com

\begin{abstract}
This article computes the number of fuzzy subgroups of symmetric group $ S_{5 }$  . The\emph{ Inclusion-Exclusion principle} is used to  determine the number of distinct fuzzy subgroups of \textbf{ symmetric group} $ S_{5} $  . Some inequalities satisfied by this number  are also established for $ n \geq 5 $


\end{abstract}

\section{Introduction}
The concept of fuzzy sets was first introduced by Zadeh in 1965. The study of fuzzy algebraic structures was started with the introduction of the concept of fuzzy subgroups by Rosenfeld in 1971. Since the first paper by Rosenfeld, researchers have sought to characterize the fuzzy subgroups of various groups.\\ One of the most important problem  of fuzzy  theory is to classify the fuzzy subgroups of a finite groups. This topic has enjoyed a rapid development in the last few years . In our case the corresponding equivalence classes of fuzzy subgroups are closely connected to the chains of subgroups in $G$. As a guiding principle in determining the number of these classes, an essential role in solving our counting problem is computed by the Inclusion-Exclusion Principle.\\
 Tarnauceanu [3] have also computed the number of fuzzy subgroups of symmetric group $S_{4}$ by the inclusion -Exclusion Principle. Symmetric groups are probably the most important in group theory, because any finite group can be embedded in such a group. They also have remarkable applications in graph theory, in enumerative combinatorics, as well as in many branches of informatics.
\title
{\bf {\begin{center}2. \;  PRELIMINARIES \end{center}}}
Suppose that $(G, \cdot, e)$ is a group with identity $e$.  Let $S(G)$ denote the collection of all fuzzy subsets of $G$.  An element $\lambda \in S(G)$ is said to be a fuzzy subgroup of $G$ if the following two conditions are sat.
\begin{enumerate}
	\item [(i)] $\lambda(ab) \geq \in \{\lambda(a), \lambda(b)\},\;\;\forall\; a,b\in G$;
	\item [(ii)] $\lambda(a^{-1} \geq \lambda(a)$ for any $a \in G$.
\end{enumerate}
And, since $(a^{-1})^{-1}=a$, we have that $\lambda(a^{-1}) = \lambda(a)$, for any $a\in G$.\\
Also, by this notation and definition, $\lambda(e) = \sup\lambda(G)$.\\
{\bf Theorem :}The set $FL(G)$ possessing all fuzzy subgroups of $G$ forms a lattice under the usual ordering of fuzzy set inclusion.  This is called the fuzzy subgroup lattice of $G$.
This theorem gives a link between  $FL(G)$ and $L(G)$, the classical subgroup lattice of $G$

We define the level subset:
$$\lambda G_{\beta} = \{a\in G/\lambda(a) \geq \beta\}\;\;\mbox{for each $\beta\in [0,1]$}$$
The fuzzy subgroups of a finite group $G$ are thus, characterized based on these level subsets.  In the sequel, $\lambda$ is a fuzzy subgroup of $G$ if and only if its level subsets are subgroups in $G$.
.

Moreover, some natural relations on $S(G)$ can also be used in the process of classifying the fuzzy subgroups of a finite  $p$-group $G$ .  One of them is defined by: $\lambda\sim \gamma$ iff $(\lambda(a)>\lambda(b) \Longleftrightarrow v(a)>v(b),\;\;\forall\;a,b\in G)$.  Also, two fuzzy subgroups $\lambda, \gamma$ of $G$ and said to be distinct if $\lambda \times v$.

As a result of this development, let $G$ be a finite group and suppose that $\lambda :G \longrightarrow [0, 1]$ is a fuzzy subgroup of $G$.  Put $\lambda(G) = \{\beta_1,\beta_2,\dots, \beta_k\}$ with the assumption that $\beta_1 < \beta_2 > \cdots > \beta_k$.  Then, ends in $G$ is determined by $\lambda$.
$$\lambda G_{\beta_1} \subset \lambda G_{\beta_2} \subset \cdots \subset \lambda G_{\beta_k} = G $$
Also, we have that:
$$\lambda(a) = \beta_t \Longleftrightarrow t = \max\{r/a \in \lambda G_{\beta_r}\}\Longleftrightarrow a\in \lambda G_{\beta_t} \backslash \lambda G_{\beta_{t-1}},\;$$
for any $a\in G$ and $t = 1,\dots, k$, where by  convention, set $\lambda G_{\beta_0}=\phi$.
\ \\ \ \\
 Hence there exits a one-to-one correspondence between the collection of the equivalence classes of fuzzy subgroups of  $G$ and the collection of chains of subgroups of $G$ which end in $G$. So,the problem of counting all distinct fuzzy subgroups of $G$ can be translated into a combinatorial problem on the subgroup lattice $ L(G) $  of  $G$ .\\  

In order to compute  the number of all distinct fuzzy subgroups of a finite $G$ which is denoted by $ h( G )$, we shall apply the inclusion-Exclusion Principle.( see[3])
$$h(G) = 2\left(\sum^t_{r=1}h(M_r) - \sum_{1\leq r_1<r_2\leq t}h(M_{r_1}\cap M_{r_2})\right.
\left.+\cdots+(-1)^{t-1}h\left(\bigcap^t_{r=1}M_r\right)\right)\hfill (\#)$$
One way to investigate the structure of a finite group is to study its maximal subgroups. In nontrivial finite groups , maximal subgroups will always exist because  the subgroups form a partially ordered set under inclusion. Since the set of subgroups is finite , this partially ordered set will have a maximal element.\\  

\section{Main Results}

 In particular , we have $h(D_{4}) = 8 $, $h(D_{6}) = 10$ , $h(D_{8}) = 32$ , $h(D_{10}) = 14 $, $h(D_{12}) = 68$ and  $h(D_{20}) = 48$\\ $h(C_2\times C_2)\cong h(D_{4}) = 8 $\\ $h(A_{3})\cong h(C_{3}) = 2 $\\
{\bf   Theorem[1] :} The number $h(A_{5})$ of all distinct fuzzy subgroups of the alternating group $A_{5}$ is 408\\
It  follows$(\#)$ that
\[ c_r = (-1)^{r-1}\sum_{0\leq i_1 < i_2\dots < i_r\leq 22} h(Mi_1\cap Mi_2\cap\dots\cap M_{r})\]
 
\[ Mi_1 \cap Mi_2 \cap \dots \cap Mi_r = \{ e \},\mbox{ for all }  r \geq 8 \mbox{ and  all } 0 \leq i_1< i_2<\dots< i_r\leq 22  \]
We have 
$$\left.\begin{array}{l}
\disp c_1=(h(A_{5})) + 5h(S_{4}) + 6h(D_{20}) + 10h(D_{12}) = 2536\\
\disp c_2= - (5h(A_{4}) +  40h(S_{3}) + 6h(D_{10}) + 45h(C_2\times C_2) + 45h(C_4) + 90h(C_2) ) = - 1324 \\
 \disp c_3 = \left( {\begin{array}{*{20} c }  22 \\ 3 \\  \end{array} } \right ) - 630 + 10h(S_{3})+ 15h(C_{4}) + 15h(C_2\times C_2) + 30h(C_{3}) + 560h(C_{2}) = 2370 \\
\disp c_4 = - \left( {\begin{array}{*{20} c }  22 \\ 4 \\  \end{array} } \right )- 586 + 576h(C_{2})+ 10h(C_3) = - 7901 
\end{array}\right.$$
$$\left.\begin{array}{ll}
c_5 = \left( {\begin{array}{*{20} c }  22 \\ 5 \\  \end{array} } \right ) - 300 + 300h(C_{2}) = 26634 &
\disp c_6 =- \left( {\begin{array}{*{20} c }  22 \\ 6 \\  \end{array} } \right ) - 85 + 85h(C_{2}) = -74698 \\
\disp c_7 = \left( {\begin{array}{*{20} c }  22 \\ 7 \\  \end{array} } \right ) - 10 + 10h(C_{2}) = 170555 &
\disp c_8 = - \left( {\begin{array}{*{20} c }  22 \\ 8 \\  \end{array} } \right ) = - 319770 \\
\disp c_9 = \left( {\begin{array}{*{20} c }  22 \\ 9 \\  \end{array} } \right ) = 497420 &
\disp c_{10} = - \left( {\begin{array}{*{20} c }  22 \\ 10 \\  \end{array} } \right ) =  - 646646 \\
\disp c_{11} = \left( {\begin{array}{*{20} c }  22 \\ 11 \\  \end{array} } \right ) = 705432 &
\disp c_{12} = - \left( {\begin{array}{*{20} c }  22 \\ 12 \\  \end{array} } \right ) = - 646646 \\
\disp c_{13} = \left( {\begin{array}{*{20} c }  22 \\  13 \\  \end{array} } \right ) = 497420 &
\disp c_{14} = - \left( {\begin{array}{*{20} c }  22 \\ 14 \\  \end{array} } \right ) = -319770 \\
\disp c_{15} = \left( {\begin{array}{*{20} c }  22 \\ 15 \\  \end{array} } \right ) = 170544 &
\disp c_{16} = - \left( {\begin{array}{*{20} c }  22 \\ 16 \\  \end{array} } \right ) = - 74613 \\
\disp c_{17} = \left( {\begin{array}{*{20} c }  22 \\ 17 \\  \end{array} } \right )  = 26334 &
\disp c_{18} = - \left( {\begin{array}{*{20} c }  22 \\ 18 \\  \end{array} } \right ) = -7315 \\
\disp c_{19} = \left( {\begin{array}{*{20} c }  22 \\ 19 \\  \end{array} } \right ) = 1540 &
\disp c_{20} = - \left( {\begin{array}{*{20} c }  22 \\ 20 \\  \end{array} } \right ) = -231 \\
\disp c_{21} = \left( {\begin{array}{*{20} c }  22 \\ 21 \\  \end{array} } \right ) = 22 &
\disp c_{22} = - \left( {\begin{array}{*{20} c }  22 \\ 22 \\  \end{array} } \right ) = -1 \\ 
\disp h(S_{5}) = 2 \sum_{r = 1}^{22} c_{r} = 3784 &\\
\end{array}\right.$$
{\bf   Theorem :} The  number $h(S_{5})$ of all distinct fuzzy subgroups of the symmetric group $S_{5}$ is 3784.\\ 
{\bf   Theorem[3] :} For $ n\geq 5$, the number $h(S_{n})$ of all distinct fuzzy subgroups of the symmetric group $ S_{n}$ satisfies the following inequality;

$$h(S_{n})\geq 2\left(\sum_{r=0}^{n}(- 1)^{r}\left( {\begin{array}{*{20} c }  n \\ r \\  \end{array} } \right )h(A_{n-r}) + \sum_{r=0}^{n - 1}\left( {\begin{array}{*{20} c }  n\\ r +1 \\  \end{array} } \right )h(S_{n - r - 1})  \right)$$

 Such a bound can be also inferred for the number of all distinct fuzzy subgroups of $(S_{5})$, we have  $h(S_{5})\geq 1940 +2h(A_{5})$ \\ It implies that,
$h(S_{5}) > 1948$

\newpage

\providecommand*{\bibname}{}
\renewcommand*{\bibname}{References}

\end{document}